\documentclass[11pt]{amsart}
\usepackage{amssymb,amsmath,txfonts,mathrsfs,mathtools,titletoc, tikz, stmaryrd}
\usepackage[alphabetic]{amsrefs}
\usepackage{color}

\newtheorem{theorem}{Theorem}[section]
\newtheorem{prop}[theorem]{Proposition}

\newtheorem{remark}[theorem]{Remark}

\numberwithin{equation}{section}

\def\pf{{\it Proof:}~}
\def\R{{\bf\mathbb R} }
\begin{document}

\title[The sharp lower bound of the first Dirichlet eigenvalue]{The sharp lower bound of the first Dirichlet eigenvalue for geodesic balls}
\author{Haibin Wang, Guoyi Xu, Jie Zhou}
\address{Haibin Wang\\ Department of Mathematical Sciences\\Tsinghua University, Beijing\\P. R. China}
\email{wanghb20@mails.tsinghua.edu.cn}
\address{Guoyi Xu\\ Department of Mathematical Sciences\\Tsinghua University, Beijing\\P. R. China}
\email{guoyixu@tsinghua.edu.cn}
\address{Jie Zhou\\ Department of Mathematical Sciences\\Tsinghua University; Academy for Multidisciplinary Studies\\ Capital Normal University,  Beijing\\P. R. China}
\email{zhoujie2014@mails.ucas.ac.cn}
\date{\today}
\date{\today}

\begin{abstract}
On complete noncompact Riemannian manifolds with non-negative Ricci curvature, Li-Schoen proved the uniform Poincar\'{e} inequality for any geodesic ball. In this note, we obtain the sharp lower bound of the first Dirichlet eigenvalue of such geodesic balls, which implies the sharp Poincar\'{e} inequality for geodesic balls.
\\[3mm]
\end{abstract}
\thanks{The second author was partially supported by NSFC 11771230 , NSFC 12026409 and Beijing Natural Science Foundation Z190003}

\maketitle

\section{Introduction}

In $1984$, Li and Schoen \cite{LS} proved the uniform Poincar\'e inequality for geodesic balls of complete noncompact Riemannian manifold with $Rc\geq 0$, their method is the integration estimate. Their Poincar\'{e} inequality is not sharp. In this note, we prove the sharp Poincar\'{e} inequality by obtaining the the sharp lower bound of the first Dirichlet eigenvalue of geodesic balls as follows.
\begin{theorem}\label{thm sharp Poincare ineq-intro}
{Suppose that $M^n$ is a complete non-compact Riemannian manifold with $Rc\geq 0$, for any $p\in M^n$, we have $\lambda_1(B_1(p))> \frac{\pi^2}{16}$, where $B_1(p)\subseteq M^n$ is the unit geodesic ball centered at $p$. Furthermore, this lower bound is sharp.
}
\end{theorem}

\begin{remark}\label{rem sharp bound}
{The above lower bound is sharp because we find a sequence of Riemannian manifolds $(M^n_i, g_i)$ satisfying the above assumption and $p_i\in M_i^n$, such that $\displaystyle \lim_{i\rightarrow\infty}\lambda_1(B_1(p_i))= \frac{\pi^2}{16}$ (see Figure \ref{figure: lemma2.1}).
}
\end{remark}

\begin{figure}
\begin{center}
\includegraphics{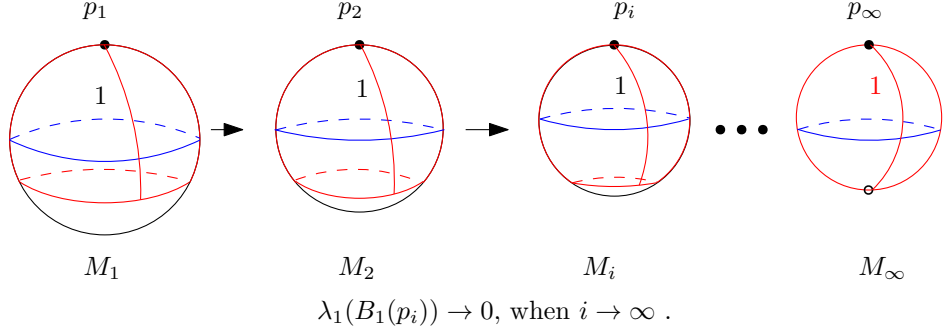}
\caption{Eigenvalues of geodesic balls in compact manifolds}
\label{figure: bigcap}
\end{center}
\end{figure}

The above sharp lower bound holds only for geodesic balls on complete non-compact Riemannian manifolds, because we know the following example from \cite{CF78}. Consider $M_i= \mathbb{S}^n(\frac{1+ 2^{-i}}{\pi}) \subseteq \R^{n+ 1}$, which is the boundary of the ball with radius $\frac{1+ 2^{-i}}{\pi}$ in $\mathbb{R}^{n+ 1}$; then $\displaystyle \lim_{i\rightarrow\infty}\lambda_1(B_{1}(p_i)) = 0$, where $p_i$ is the north pole of $M_i$ and $B_{1}(p_i)$ is the unit geodesic ball centered at $p_i$ (see Figure \ref{figure: bigcap}).

One well-known tool to get the sharp estimate of eigenvalues on manifolds is, the gradient estimate of eigenfunction introduced by Li-Yau \cite{LY}. By establishing the sharp gradient estimate of Neumann eigenfunction, Zhong-Yang \cite{ZY} proved that $\mu_2(M^n)\geq \frac{\pi^2}{4}$ for compact Riemannian manifold $M^n$ with $Rc\geq 0$ and $\mathrm{diam}(M^n)= 2$, where $\mu_2(M^n)$ is the second Neumann eigenvalue of $M^n$ and $\mu_1(M^n)= 0$. Furthermore, Hang-Wang \cite{HW} showed the equality holds if and only if $M^n$ is isometric to $\mathbb{S}^1$. Later, Yang \cite{Yang} applied the gradient estimate method to obtain the lower bound of the first Dirichlet eigenvalue for domain with mean convex boundary.

For any domain $\Omega\subseteq (M^n, g)$ with $Rc\geq 0$ and maximal volume growth (i.e. $\displaystyle \mathrm{AVR}(M^n)= \lim_{r\rightarrow \infty} \frac{\mathrm{Vol}(B_r(p))}{\omega_n r^n}> 0$, where $\omega_n$ is the volume of unit ball in $\mathbb{R}^n$); using Brendle's sharp isoperimetric inequality (see \cite{Brendle}), recently Krist\'{a}ly \cite{Kristaly} proved the sharp lower bound of the first Dirichlet eigenvalue as follows:
\begin{align}
\displaystyle \lambda_1(\Omega)\geq j_{\frac{n}{2}- 1}^2(\omega_n\mathrm{AVR}(g))^{\frac{2}{n}}\mathrm{Vol}(\Omega)^{-\frac{2}{n}}, \nonumber
\end{align}
where $j_v$ is the first positive root of the Bessel function $J_v$ of the first kind with degree $v\in \mathbb{R}$. Especially, the equality holds if and only if $(M^n, g)$ is $\mathbb{R}^n$.

Our idea is applying the maximum principle on the suitable function related to eigenfunction, and we obtain the $C^0$ estimate of Dirichlet eigenfunction directly without involving the gradient. This argument is partly inspired by the modulus continuity estimate of eigenfunctions in \cite{Ni}, and it also has close relationship with the original parabolic method of \cite{AC} (also see \cite{AC-1}). As we know, the modulus continuity estimate of eigenfunctions can not yield the $C^0$ estimate of eigenfunction in our case, and is not enough to yield $\displaystyle \lambda_1(B_1(p))> \frac{\pi^2}{16}$. This strict sharp inequality is novel, comparing the facts that the sharp lower bound of $\mu_2(M^n)$ can be realized by $M^n= \mathbb{S}^1$ in \cite{ZY} and the sharp lower bound of $\lambda_1(\Omega)$ can be realized by $\Omega= B_1^n\subseteq \mathbb{R}^n$ in \cite{Kristaly}.

\section{The sharp lower bound of the first Dirichlet eigenvalue}

On complete non-compact Riemannian manifold $M^n$, let $\gamma$ be a ray starting at $p$, define the Busemann function with respect to $\gamma$ as $\displaystyle b_\gamma(x)= \lim_{t\rightarrow\infty} \Big[t- d(x, \gamma(t)\Big]$. When the context is clear, we also use $b$ instead of $b_\gamma$ for simplicity.

\begin{prop}\label{prop region between b level sets}
{Let $M^n$ be a complete noncompact Riemannian manifold with $Rc\geq 0$, for any ray $\gamma\subseteq M^n$ and region $\Omega\subseteq \subseteq b_\gamma^{-1}[a, a+ D]$ where $a\in \mathbb{R}, D> 0$, we have $\displaystyle \lambda_1(\Omega)\geq \frac{\pi^2}{4D^2}$.
}
\end{prop}

\begin{remark}\label{rem C0 est of Diri eigenfunc}
{By establishing the $C^0$-estimate of the first Dirichlet eigenfunction, we get the sharp estimate of $\lambda_1$. One novel thing about our $C^0$-estimate (\ref{object ineq}) is that there is no boundary assumption, although we use the noncompact property of the whole manifold $M^n$.
}
\end{remark}

\pf
{For simplicity, we use $b$ instead of $b_\gamma$ in the rest of the argument. Assume $\Delta u= -\lambda_1 u$ for some constant $\lambda_1> 0$ with $\displaystyle u\in C_0^\infty(\Omega), u(x)\geq 0, \max_{x\in \Omega}u(x)= 1$. Let $\alpha= a+ D$, firstly we show that
\begin{align}
\sin^{-1}u\leq \sqrt{\lambda_1}(\alpha- b(x)). \label{object ineq}
\end{align}

To prove (\ref{object ineq}), we only need to show for any $\delta\in (0, 1)$ such that
\begin{align}
\sin^{-1} ((1- \delta)u)\leq \sqrt{\lambda_1}(\alpha- b(x)). \label{object ineq-1}
\end{align}
Then (\ref{object ineq}) follows from letting $\delta\rightarrow 0$ in (\ref{object ineq-1}).

By contradiction. If (\ref{object ineq-1}) does not hold for some $\delta\in (0, 1)$, then there is some $\epsilon> 0$ such that for $\displaystyle \phi(x)= \sin^{-1}((1- \delta)u(x))- \sqrt{\lambda_1+ \epsilon}(\alpha- b(x))$, we have $\displaystyle \max_{x\in \Omega}\phi(x)> 0$. Note $\displaystyle\max_{x\in \partial \Omega}\phi(x)= -\sqrt{\lambda_1+ \epsilon}(\alpha- b(x))\leq 0$, so there is $x_0\in \Omega$ such that $\displaystyle \phi(x_0)= \max_{x\in \Omega}\phi(x)> 0$.

Now we have $\nabla \phi(x_0)=0$, which implies $\displaystyle \frac{(1- \delta)|\nabla u|}{\sqrt{1- ((1- \delta)u)^2}}= \sqrt{\lambda_1+ \epsilon}$. Also note $\Delta b\geq 0$ by $Rc\geq 0$ and the Laplace Comparison Theorem. Now we have
\begin{align}
0&\geq \Delta\phi(x_0)= \frac{(1- \delta)u}{\sqrt{1- ((1- \delta)u)^2}}(-\lambda_1+ \frac{(1- \delta)^2|\nabla u|^2}{1- ((1- \delta)u)^2})+ \sqrt{\lambda_1+ \epsilon}\Delta b \nonumber \\
&\geq \frac{\epsilon(1- \delta)u}{\sqrt{1- ((1- \delta)u)^2}}> 0. \nonumber
\end{align}
This is the contradiction, then (\ref{object ineq-1}) is proved.

From $b(x)\geq a$, we have $0\leq \alpha- b(x)\leq D$ for any $x\in \overline{\Omega}$. Assume $u(x_1)= 1$ for some $x_1\in \Omega$, we have $\displaystyle \lambda_1\geq (\frac{\sin^{-1}u(x_1)}{\alpha- b(x_1)})^2\geq \frac{\pi^2}{4D^2}$.
}
\qed

Note $B_1(p)\subseteq b_\gamma^{-1}[-1, 1]$ for some ray starting from $p$,  we have the following sharp lower bound of the first Dirichlet eigenvalue.
\begin{theorem}\label{thm Diri eigenvalue for set with diam 2}
{Let $M^n$ be a complete noncompact Riemannian manifold with $Rc\geq 0$, for any region $\Omega\subseteq M^n$ with $\mathrm{diam}(\Omega)= D< \infty$,  we have $\displaystyle \lambda_1(\Omega)\geq \frac{\pi^2}{4D^2}$. Furthermore $\lambda_1(B_1(p))> \frac{\pi^2}{16}$ and this inequality is sharp.
}
\end{theorem}

\pf
{\textbf{Step (1)}. For $p\in \Omega$, take a ray $\gamma$ starting from $p$, let the Busemann function with respect to $\gamma$ be $b(x)$. From $\mathrm{diam}(\Omega)= D$ and $|\nabla b|= 1$, we have $\displaystyle \max_{x\in \Omega}b(x)- \min_{x\in \Omega}b(x)\leq D$, hence $\Omega\subseteq b^{-1}[a, a+ D]$ for some $a$. Then $\displaystyle \lambda_1(\Omega)\geq \frac{\pi^2}{4D^2}$ follows from Proposition \ref{prop region between b level sets}.

Assume $\Delta u= -\lambda_1 u$ for some constant $\lambda_1> 0$ with $\displaystyle u\in C_0^\infty(B_1(p)), u(x)\geq 0, \max_{x\in B_1(p)}u(x)= 1$.
To prove $\lambda_1> \frac{\pi^2}{16}$, from the above, we only need to show that $\lambda_1\neq \frac{\pi^2}{16}$. If $\lambda_1= \frac{\pi^2}{16}$, assume $u(x_1)= 1$ for some $x_1\in B_1(p)$. Then from (\ref{object ineq}), we have
\begin{align}
\frac{\pi}{4}(1- b(x_1))\geq \sin^{-1}u(x_1)= \frac{\pi}{2}. \nonumber
\end{align}
This implies $b(x_1)= -1$. On the other hand, from the definition of $b(x)$, we know $b(x_1)\geq -d(x_1, p)> -1$, it is the contradiction.

\textbf{Step (2)}. Next we show that $\displaystyle \lambda_1(B_1(p))> \frac{\pi^2}{16}$ is sharp. Consider $(\mathbb{R}^n, g_i)$ with $g_i= dr^2+ f_i^2(r)d\mathbb{S}^{n- 1}$, where $\epsilon_i= 2^{-i}$ and
\begin{equation}\nonumber
f_i(r)=\left\{
\begin{array}{rl}
&\epsilon_i- \epsilon_i e^{\epsilon_i^{-1}} e^{\frac{\epsilon_i^{-1}}{\epsilon_i^{-1} r- 1}} \quad \quad \quad \quad \quad 0\leq r< \epsilon_i ,\\
&\epsilon_i \quad \quad \quad \quad \quad \quad \quad \quad  r\geq \epsilon_i . \nonumber
\end{array} \right.
\end{equation}
Then it is easy to check that $Rc(g_i)\geq 0$. Assume the origin point of $(\mathbb{R}^n, g_i)$ is $q_i$.

\begin{figure}
\begin{center}
\includegraphics{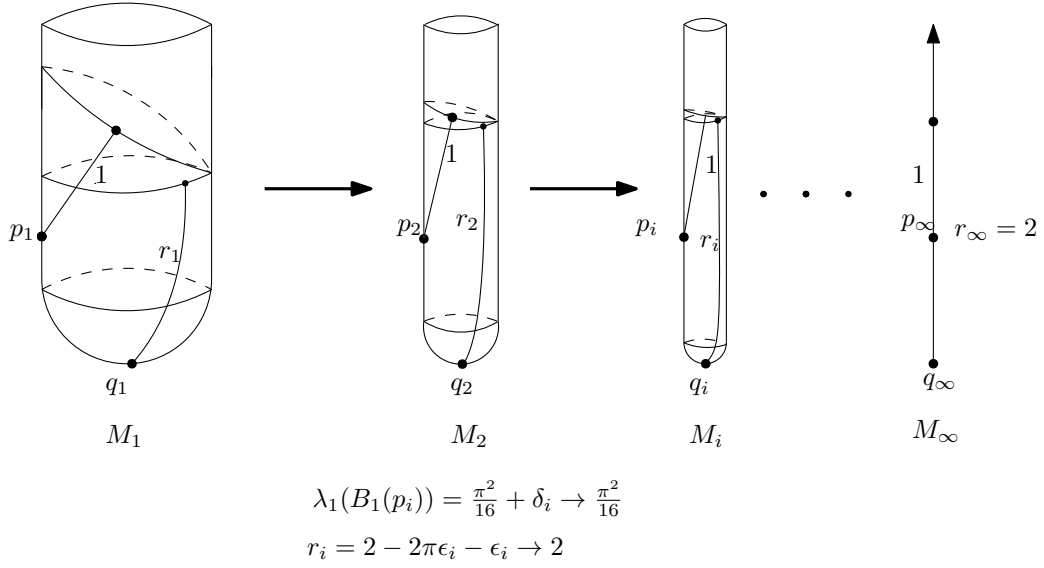}
\caption{Collapsing to the sharp lower bound}
\label{figure: lemma2.1}
\end{center}
\end{figure}

Take $p_i=\exp_{q_i}((1-\epsilon_i)\theta_0)$ for some $\theta_0\in S^{n- 1}$, we claim
$$B_{2-2\pi\epsilon_i-\epsilon_i}(q_i)\subset B_1(p_i).$$
For any $x\in B_{2-2\pi\epsilon_i-\epsilon_i }(q_i)$, we know $x=\exp_{q_i}(r_x\theta_x)$ for some $r_x\le 2-2\pi\epsilon_i-\epsilon_i$ and $\theta_x\in S^1$.  Take $x_1= \exp_{q_i}((1-\epsilon_i)\theta_x)$, by $f_i(r)\le \epsilon_i$, we know
$$d(x,p_i)\le d(x,x_1)+d(x_1,p_i)\le 1-2\pi\epsilon_i+\pi \epsilon_i<1.$$
This implies $B_{2-2\pi\epsilon_i-\epsilon_i}(q_i)\subset B_1(p_i)$. So, by the monotonicity of Dirichlet eigenvalue with respect to the domain, we know $\lambda_1(B_1(p_i))\le \lambda_1(B_{2-2\pi\epsilon_i-\epsilon_i}(q_i)).$

On the other hand, for $r_i=2-2\pi \epsilon_i-\epsilon_i$,  by taking $\tau_i=\frac{\pi}{2r_i}$ and $\varphi_i(x)=\cos{\tau_i d(x,q_i)}\in W_0^{1,2}(B_{r_i}(q_i))$. We know
\begin{align*}
\lambda_1(B_{2-2\pi\epsilon_i-\epsilon_i}(q_i))
\le \frac{\int_{B_{r_i}(q_i)}|\nabla \varphi_i|^2d\mu_{g_i}}{\int_{B_{r_i}(q_i)}\varphi_i^2d\mu_{g_i}}=\frac{\tau_i^2\int_0^{\frac{\pi}{2}}\sin^2(s)f_i^{n- 1}(\frac{s}{\tau_i})ds}{\int_0^{\frac{\pi}{2}}\cos^2(s)f_i^{n- 1}(\frac{s}{\tau_i})ds}\le \frac{\epsilon_i^{n- 1} \int_0^{\frac{\pi}{2}}\sin^2(s)ds}{\epsilon_i^{n- 1}\int_{\tau_i \epsilon_i}^{\frac{\pi}{2}}\cos^2(s)ds}\tau_i^2.
\end{align*}
Thus $\displaystyle \lim_{i\to \infty}\lambda_1(B_1(p_i))\le \lim_{i\to \infty}\frac{\int_0^{\frac{\pi}{2}}\sin^2(s)ds}{\int_{\tau_i \epsilon_i}^{\frac{\pi}{2}}\cos^2(s)ds}\tau_i^2=\frac{\pi^2}{16}$.
}
\qed

  Similar to Sakai's conjecture for Zhong-Yang's sharp eigenvalue estimate (see \cite{Sakai}) and Hang-Wang's rigidity result \cite{HW}, we may ask whether the unit geodesic ball in complete noncompact manifolds with nonnegative Ricci curvature is close to an interval $[0, 2]$ in the Gromov-Hausdorff sense if the first eigenvalue is nearly $\frac{\pi^2}{16}$. From (\ref{object ineq}), we in fact have $\displaystyle \lim_{\lambda_1(B_1(p_i))\rightarrow \frac{\pi^2}{16}}\mathrm{diam}(B_1(p_i))= 2$. One natural question is whether the width of $B_1(p)$ is nearly $0$ when $\lambda_1(B_1(p))$ is nearly $\frac{\pi^2}{16}$.

  \section*{Acknowledgments}
We thank one anonymous referee for the helpful comment on the earlier version of this note.

\end{document}